# THE DISORDER PROBLEM FOR COMPOUND POISSON PROCESSES WITH EXPONENTIAL JUMPS

By Pavel V. Gapeev

*Russian Academy of Sciences*

The problem of disorder seeks to determine a stopping time which is as close as possible to the unknown time of "disorder" when the observed process changes its probability characteristics. We give a partial answer to this question for some special cases of Lévy processes and present a complete solution of the Bayesian and variational problem for a compound Poisson process with exponential jumps. The method of proof is based on reducing the Bayesian problem to an integro-differential free-boundary problem where, in some cases, the smooth-fit principle breaks down and is replaced by the principle of continuous fit.

**1. Introduction.** Assume that at time $t = 0$ we begin to observe a continuously updated process $X = (X_t)_{t \geq 0}$ whose probability characteristics change at some unknown time $\theta$, called the *time of disorder*, which cannot be observed directly. Throughout this paper the random time $\theta$ can take the value 0 with probability $\pi$; under the condition that $\theta > 0$, it is exponentially distributed with parameter $\lambda > 0$. The disorder problem or the problem of quickest disorder detection is to decide by observing the process $X$ the time instant at which we should give an alarm to indicate the occurrence of disorder. This time instant should be as close as possible to the time $\theta$ in the sense that both the probability of false alarm and the expectation of the time interval between the occurrence of disorder and the alarm (when the latter is given correctly) should be minimal.

The problem of detecting a change in drift of a Wiener process was formulated and solved by Shiryaev [12, 13, 14, 15] (see also [16] and [17], Chapter IV and page 208, for historical notes and references). Some particular cases of the problem of detecting a change in the intensity of a Poisson process were considered by Gal'chuk and Rozovskii [6] and by Davis [4]. Peskir and Shiryaev [10] presented a complete solution of the disorder problem for a Poisson process in the Bayesian formulation. The main aim of this paper is to find an explicit expression of the optimal stopping boundary for the a posteriori probability process in some special cases of the problem for Lévy processes and to present a complete solution to the problem for a compound Poisson process that has exponentially distributed jumps. Actually, we give the next example of process for which the quickest disorder detection problem can be solved in an explicit form. Such processes are used, for example, in several models of stochastic finance and insurance (see, e.g., [18]). For some other optimal stopping problems for such processes see, for example, [9].

The paper is organized as follows. In Section 2 we give a formulation of the Bayesian and variational problem of quickest disorder detection for Lévy processes. In Section 3 by the examination of the sample-path behavior of the a posteriori probability process, we find an optimal stopping boundary in some particular cases of the Bayesian problem. In Section 4 by means of solving the corresponding integro-differential free-boundary problem, we derive a complete solution of the Bayesian problem for a compound Poisson process with exponential jumps, where we can observe the breakdown of the smooth-fit principle and its replacement by the principle of continuous fit. These effects can be explained both by the examination of the sample-path properties of the a posteriori probability process and by the existence of a singularity point of the integro-differential equation. Note that in models based on jump processes the situations when the continuous fit replaces the smooth fit were earlier observed, for example, in bandit problems (see, e.g., [2] for references). In Section 5, passing from the derived solution of the Bayesian problem, we find an explicit expression for the optimal stopping boundary in the corresponding variational problem.

We note here that the problem of quickest detection admits different formulations and appears in on-line quality control, radar location, seismology and so forth (see, e.g., [3, 8]).

**2. Formulation of the Bayesian and variational problem.** For a precise probabilistic formulation of the quickest disorder detection problem for Lévy processes (see [17], Chapter IV, for the Wiener process case), let us suppose that on some measurable space $(\Omega, \mathcal{F})$ equipped with a family of probability measures $(P^s)_{s \geq 0}$ there exists a nonnegative random variable $\theta$ such that $P^s[\theta = s] = 1$ for all $s \geq 0$. It is assumed that we observe a continuously updated process $X = (X_t)_{t \geq 0}$ with $X_0 = 0$ and having, under the measure $P^s$, the triplet

$$(2.1) \quad ((t \wedge s)b_0 + ((t-s) \vee 0)b_1, 0, dt\,[I_{\{t<s\}}\nu_0(dx) + I_{\{t \geq s\}}\nu_1(dx)])$$

with respect to the function $h(x) = x$, $x \in \mathbb{R}$, for all $t, s \geq 0$, where $\nu_i(dx)$ is a Lévy measure on $\mathbb{R}$ such that $\nu_i(\{0\}) = 0$ and $\int (x^2 \wedge 1)\nu_i(dx) < \infty$ for $i = 0, 1$ (see, e.g., [7], Chapter II.4, or [11], Chapter II.8). Here $\theta$ and $X$ are assumed to be stochastically independent under $P^s$ for all $s \geq 0$. Let us fix $\lambda > 0$ and define the measures $P_\pi = \pi P^0 + (1 - \pi) \int_0^\infty \lambda e^{-\lambda s} P^s \, ds$ for all $\pi \in [0, 1]$, so that we have $P_\pi[\theta = 0] = \pi$ and $P_\pi[\theta > t | \theta > 0] = e^{-\lambda t}$ for all $t \geq 0$.

Let $\tau$ be a stopping time with respect to the filtration $\mathbf{F}^X = (\mathcal{F}_t^X)_{t \geq 0}$, where $\mathcal{F}_t^X = \sigma\{X_s | 0 \leq s \leq t\}$. We interpret $\tau$ as the time at which the alarm is sounded to signal the change in distribution of the observed process $X$. The *Bayesian disorder problem* is to minimize the risk function

$$(2.2) \qquad V(\pi) = \inf_\tau \{P_\pi[\tau < \theta] + cE_\pi[\tau - \theta]^+\},$$

where the infimum is taken over all $\mathbf{F}^X$ stopping times $\tau$, and to find an optimal stopping time $\tau_*$ at which the infimum in (2.2) is attained. Here $P_\pi[\tau < \theta]$ is the probability of false alarm, $E_\pi[\tau - \theta]^+$ is the average delay in detecting disorder correctly and $c > 0$ is some constant.

It is easily shown (see [17], pages 195–197) that the value function $V(\pi)$ can be expressed in terms of the a posteriori probability process $(\pi_t)$, where $\pi_t = P_\pi[\theta \leq t | \mathcal{F}_t^X]$ for all $t \geq 0$ and $P_\pi[\pi_0 = \pi] = 1$. Namely, we have

$$(2.3) \qquad V(\pi) = \inf_\tau E_\pi \left[ 1 - \pi_\tau + c \int_0^\tau \pi_t \, dt \right].$$

Moreover, it is easily verified (see [17], page 204) that the infimum in (2.3) is actually taken over the class $\mathcal{M}(\pi)$ of stopping times $\tau$ such that $E_\pi[\tau] < \infty$.

To give the corresponding variational or fixed false-alarm probability formulation, let the number $\pi \in [0, 1)$ be fixed and let $\mathcal{M}(\pi, \alpha)$ denote the class of stopping times $\tau$ that satisfy

$$(2.4) \qquad P_\pi[\tau < \theta] \leq \alpha,$$

where $\alpha$ is a given constant from the interval $[0, 1)$. The *variational disorder problem* is to find in the class $\mathcal{M}(\pi, \alpha)$ a stopping time $\hat{\tau}$ such that

$$(2.5) \qquad E_\pi[\hat{\tau} - \theta]^+ \leq E_\pi[\tau - \theta]^+$$

for any other stopping time $\tau$ from $\mathcal{M}(\pi, \alpha)$.

**3. Preliminary results and examples.** Suppose that the filtration $\mathbf{F}^X$ is right-continuous and the conditions

$$(3.1) \qquad \int |x|\nu_i(dx) < \infty \qquad (i = 0, 1),$$

$$(3.2) \qquad b_1 = b_0 + \int x\nu_1(dx) - \int x\nu_0(dx),$$

$$(3.3) \qquad \int (\sqrt{Y(x)} - 1)^2 \nu_0(dx) < \infty$$

are satisfied, where $Y(x) = \nu_1(dx)/\nu_0(dx)$ for all $x \in \mathbb{R}$. Then by means of Girsanov's theorem for semimartingales ([7], Theorem III.5.34) and Itô's formula ([7], Theorem I.4.57), using the schema of arguments in [17], page 202, it can be verified that the process $(\pi_t)$ solves the stochastic differential equation

$$(3.4) \quad d\pi_t = \lambda(1 - \pi_t)\, dt + \int \frac{\pi_{t-}(1 - \pi_{t-})(Y(x) - 1)}{1 + \pi_{t-}(Y(x) - 1)} (\mu^X - \nu^X)(dt, dx),$$

where $\mu^X$ is the measure of jumps of the process $X$ and its $\mathbf{F}^X$ compensator $\nu^X$ is given by $\nu^X(dt, dx) = (\pi_{t-}\nu_1(dx) + (1 - \pi_{t-})\nu_0(dx))\, dt$. From (3.4) it is easily seen that $(\pi_t)$ is a time-homogeneous (strong) Markov process under $P_\pi$ with respect to the natural filtration which clearly coincides with $\mathbf{F}^X$. The latter implies that the infimum in (2.3) can be taken over all stopping times of $(\pi_t)$ playing the role of a sufficient statistic (see, e.g., [17], Chapter II.15).

It can be also verified (see [17], pages 197 and 198, and [10]) that the value function $V(\pi)$ is decreasing and concave on $[0, 1]$, and the optimal stopping time in (2.3) is given by

$$(3.5) \quad \tau_* = \inf\{t \geq 0 | \pi_t \geq B_*\},$$

where $B_*$ is the smallest number $\pi$ from $[0, 1]$ such that $V(\pi) = 1 - \pi$.

Using the arguments from [10] we now find an explicit expression for the optimal stopping boundary $B_*$ in some particular cases of the problem.

LEMMA 3.1. *Assume in addition to* (2.1) *and* (3.1)–(3.3) *that we have*

$$(3.6) \quad \nu_1(dx) \geq \nu_0(dx) \qquad (x \in \mathbb{R}),$$

$$(3.7) \quad 0 < \int x\nu_1(dx) - \int x\nu_0(dx) \leq c + \lambda.$$

*Then in the Bayesian problem of quickest disorder detection* (2.2) + (2.3) *the stopping time $\tau_*$ from* (3.5) *is optimal with $B_* = \overline{B}$, where we set*

$$(3.8) \quad \overline{B} = \frac{\lambda}{\lambda + c}.$$

PROOF. The assumption (3.7) ensures that $\overline{B} \leq \widehat{B}$, where we set

$$(3.9) \quad \widehat{B} = \lambda \Big/ \Big( \int x\nu_1(dx) - \int x\nu_0(dx) \Big).$$

From (3.4) it is seen that if $\widehat{B} \geq 1$, then the process $(\pi_t)$ is strictly increasing, and if $\widehat{B} < 1$, then the drift rate of the continuous part of $(\pi_t)$ is positive on $[0, \widehat{B})$, negative on $(\widehat{B}, 1)$ and equal to zero at $\widehat{B}$. Thus, if $(\pi_t)$ starts in $[0, \widehat{B})$ or in $(\widehat{B}, 1)$, then under the absence of jumps, $(\pi_t)$ never reaches $\widehat{B}$, because its drift tends to zero the same time with linear order. Therefore, by virtue of the fact that under the condition (3.6) the process $(\pi_t)$ can have only positive jumps, it can leave $[0, \widehat{B})$ only by jumping and, fluctuating in $(\widehat{B}, 1)$, it cannot enter $[0, \widehat{B})$. If $(\pi_t)$ starts or ends up at $\widehat{B}$, then it is trapped there ($P_\pi$-a.s.) until the next jump of $X$ occurs.

From (3.4) it follows that the process $(\pi_t)$ admits the representation

$$\pi_t = \pi + \lambda \int_0^t (1 - \pi_{s-}) \, ds + M_t, \tag{3.10}$$

where $(M_t)$ is a martingale under $P_\pi$ with respect to $\mathbf{F}^X$. Hence, by means of the optional sampling theorem (see, e.g., [7], Theorem I.1.39), from (3.10) together with (3.4) and according to (3.1) we obtain that $E_\pi[M_\tau] = 0$ and hence

$$E_\pi \left[ 1 - \pi_\tau + c \int_0^\tau \pi_t \, dt \right] = 1 - \pi + (\lambda + c) E_\pi \int_0^\tau \left( \pi_t - \frac{\lambda}{\lambda + c} \right) dt \tag{3.11}$$

for all stopping times $\tau$ from $\mathcal{M}(\pi)$. Recalling that the process $(\pi_t)$ is monotone increasing in $[\overline{B}, \widehat{B})$ and that after entering $[\widehat{B}, 1]$ cannot leave it, from (3.11) we may therefore conclude that it is never optimal to stop $(\pi_t)$ in $[0, \overline{B})$ and that $(\pi_t)$ must be stopped instantly after passing through $\overline{B}$. □

EXAMPLE 3.2. Assume that in (2.1) we have $b_i = 1/\lambda_i$ and $\nu_i(dx) = I_{\{x > 0\}} \exp(-\lambda_i x) \, dx/x$ with $\lambda_i > 0$. Thus $X$ is a gamma process with parameter changing from $\lambda_0$ to $\lambda_1$ (see, e.g., [18], Chapter III.1). In this case the integrals in (3.1) and (3.3) are equal to $1/\lambda_i$ and $\log[(\lambda_0 + \lambda_1)^2/(4\lambda_0 \lambda_1)]$, respectively. Therefore, by Lemma 3.1 we get that if $\lambda_0 > \lambda_1 > 0$ and $\log(\lambda_0/\lambda_1) \le c + \lambda$, then the stopping time $\tau_*$ from (3.5) is optimal with $B_* = \lambda/(\lambda + c)$.

EXAMPLE 3.3. Suppose that in (2.1) we have $b_i = 1/\gamma_i$ and $\nu_i(dx) = I_{\{x > 0\}} \exp(-\gamma_i^2 x/2) \, dx/(2\pi x^3)^{1/2}$ with $\gamma_i > 0$. Thus $X$ is an inverse Gaussian process with parameter changing from $\gamma_0$ to $\gamma_1$ (see, e.g., [1]). In this case the integrals in (3.1) and (3.3) are equal to $1/\gamma_i$ and $[2(\gamma_0^2 + \gamma_1^2)]^{1/2} - \gamma_0 - \gamma_1$, respectively. Therefore, by Lemma 3.1 we conclude that if $\gamma_0 > \gamma_1 > 0$ and $\gamma_0 - \gamma_1 \le c + \lambda$, then $\tau_*$ from (3.5) is optimal with $B_* = \lambda/(\lambda + c)$.

REMARK 3.4. From (3.11) it is seen that one should not stop $(\pi_t)$ when it is in $[0, \overline{B}]$, so for $B_*$ from (3.5) we have $\overline{B} \le B_* \le 1$.

**4. Solution of the Bayesian problem for a compound Poisson process with exponential jumps.** In the rest of the paper, we assume that the process $X$ is defined by

$$X_t = \int_0^t \theta_{s-} \, dX_s^1 + \int_0^t (1 - \theta_{s-}) \, dX_s^0, \tag{4.1}$$

where $X_s^i = \sum_{j=1}^{N_s^i} \xi_j^i$ and $\theta_s = I_{\{s \geq \theta\}}$ for all $t, s \geq 0$, $N^i = (N_t^i)$ is a Poisson process with intensity $1/\lambda_i$ and $(\xi_j^i)_{j \in \mathbb{N}}$ is a sequence of independent random variables exponentially distributed with parameter $\lambda_i$ [$N^i$, $(\xi_j^i)_{j \in \mathbb{N}}$ and $\theta$ are supposed to be independent] for $i = 0, 1$. Then in (2.1) we have $b_i = 1/\lambda_i^2$ and $\nu_i(dx) = I_{\{x>0\}} \exp(-\lambda_i x) dx$, and thus $X$ is a compound Poisson process that has exponentially distributed jumps with parameter changing from $\lambda_0$ to $\lambda_1$. In this case the integrals in (3.1) and (3.3) are equal to $1/\lambda_i^2$ and $(\lambda_0 - \lambda_1)^2 / [\lambda_0 \lambda_1 (\lambda_0 + \lambda_1)]$, respectively, and (3.4) takes the form

$$
\begin{aligned}
d\pi_t = {}& \lambda(1 - \pi_t) \, dt \\
& + \int_0^\infty \frac{\pi_{t-}(1 - \pi_{t-})(\exp(-\lambda_1 x) - \exp(-\lambda_0 x))}{\pi_{t-} \exp(-\lambda_1 x) + (1 - \pi_{t-}) \exp(-\lambda_0 x)} \\
& \quad \times (\mu^X(dt, dx) - (\pi_{t-} \exp(-\lambda_1 x) \\
& \qquad + (1 - \pi_{t-}) \exp(-\lambda_0 x)) \, dt \, dx).
\end{aligned}
\tag{4.2}
$$

Standard arguments imply that in this case the infinitesimal operator $\mathbb{L}$ of the process $(\pi_t)$ acts on a function $f \in C^1([0, 1])$ according to the rule

$$
\begin{aligned}
(\mathbb{L}f)(\pi) = {}& \left( \lambda - \left( \frac{\lambda_0 - \lambda_1}{\lambda_0 \lambda_1} \right) \pi \right)(1 - \pi) f'(\pi) \\
& + \int_0^\infty \left[ f\left( \frac{\pi \exp(-\lambda_1 x)}{\pi \exp(-\lambda_1 x) + (1 - \pi) \exp(-\lambda_0 x)} \right) - f(\pi) \right] \\
& \quad \times (\pi \exp(-\lambda_1 x) + (1 - \pi) \exp(-\lambda_0 x)) \, dx
\end{aligned}
\tag{4.3}
$$

for all $\pi \in [0, 1]$. Using standard arguments based on the strong Markov property, it follows that $V(\pi)$ is $C^1$ on $(0, B_*)$. Therefore, using the results from [17], Chapter III.8, we can formulate the *integro-differential free-boundary problem* for the unknown function $V(\pi)$ from (2.3) and the unknown boundary $B_*$ from (3.5) as

$$(\mathbb{L}V)(\pi) = -c\pi \qquad (0 < \pi < B_*), \tag{4.4}$$

$$V(\pi) = 1 - \pi \qquad (B_* \leq \pi \leq 1), \tag{4.5}$$

$$V(B_*-) = 1 - B_* \qquad \text{(continuous fit)}, \tag{4.6}$$

where the condition (4.6) is satisfied by virtue of the concavity argument above. Note that the superharmonic characterization of the value function (see [5] and [17]) implies that $V(\pi)$ is the largest function that satisfies (4.4)–(4.6). Moreover, under some relationships on the parameters of the model which are specified below, the condition

$$V'(B_*) = -1 \qquad \text{(smooth fit)} \tag{4.7}$$

may be satisfied or break down. We also observe that, in this case, $\widehat{B}$ from (3.9) takes the form

$$\widehat{B} = \frac{\lambda \lambda_0 \lambda_1}{\lambda_0 - \lambda_1} \tag{4.8}$$

and turns out to be a singularity point of (4.4) when $\lambda_0 > \lambda_1$.

Using the schema of arguments in [10], we further show that the system (4.4)–(4.6) admits an explicit solution which turns out to be a solution of the initial optimal stopping problem (2.3). For this, let us consider a continuous function $f(\pi)$ that satisfies (4.4) on $(0, B)$ and (4.5) on $[B, 1]$ for some $0 < B < 1$ given and fixed.

Let us first assume that $\lambda_0 > \lambda_1$. Then it follows that the function $\tilde{f}(y) = f(\pi)$ with $\pi = e^y/(1 + e^y)$ solves the system

$$\left(\frac{\lambda'(1+e^y)}{e^y} - \frac{1}{\gamma(\gamma-1)}\right)\tilde{f}'(y) - \frac{\tilde{f}(y)[\gamma(1+e^y)-1]}{\gamma(\gamma-1)(1+e^y)}$$

$$+ \frac{e^{\gamma y}}{1+e^y}\left[\int_y^{\widetilde{B}} \frac{\tilde{f}(z)(1+e^z)}{e^{\gamma z}}\,dz + \frac{e^{-\gamma \widetilde{B}}}{\gamma}\right] \tag{4.9}$$

$$= -\frac{c(\lambda_0 - \lambda_1)e^y}{1+e^y} \qquad (y < \widetilde{B}),$$

$$\tilde{f}(y) = \frac{1}{1+e^y} \qquad (y \geq \widetilde{B}), \tag{4.10}$$

where we set $\gamma = \lambda_0/(\lambda_0 - \lambda_1) > 1$, $\lambda' = \lambda(\lambda_0 - \lambda_1) > 0$ and $\widetilde{B} = \log[B/(1-B)]$. It can be easily shown that the system (4.9) + (4.10) has a unique solution which is given by

$$\tilde{f}(y; \widetilde{B}) = \frac{1}{1+e^{\widetilde{B}}} - \int_y^{\widetilde{B}} \frac{\gamma(\gamma-1)\widetilde{F}(z, \widetilde{B})e^{\gamma z}}{\gamma(1+e^z)-1}\,dz, \tag{4.11}$$

$$\widetilde{F}(y, \widetilde{B}) = \frac{1}{\widetilde{A}(y)}\left(\widetilde{C}(y, \widetilde{B}) - \int_y^{\widetilde{B}} \frac{\widetilde{C}(z, \widetilde{B})}{\widetilde{A}(z)} \frac{\widetilde{G}(z)}{\widetilde{G}(y)}\,dz\right), \tag{4.12}$$

$$\widetilde{A}(y) = \frac{1+e^y}{e^y}\left(\frac{\lambda'\gamma(\gamma-1)(1+e^y)-e^y}{\gamma(1+e^y)-1}\right), \tag{4.13}$$

$$\widetilde{C}(y, \widetilde{B}) = \frac{e^{-(\gamma-1)\widetilde{B}}}{\gamma(\gamma-1)(1+e^{\widetilde{B}})} - \frac{c\lambda_0 e^{-(\gamma-1)y}}{\gamma}, \tag{4.14}$$

$$\widetilde{G}(y) = \begin{cases} \left|e^y - \dfrac{\widehat{B}}{1-\widehat{B}}\right|^a (1+e^y), & \text{if } \widehat{B} \neq 1, \\ \exp[-\gamma e^y](1+e^y), & \text{if } \widehat{B} = 1, \end{cases} \tag{4.15}$$

for $y \leq \widetilde{B}$, and $a = (\widehat{B} + \gamma - 1)/(1 - \widehat{B})$ if $\widehat{B} \neq 1$. Using (4.11)–(4.15) we may thus conclude that the function $f(\pi; B) = \tilde{f}(y; \widetilde{B})$ given by

$$(4.16) \quad f(\pi; B) = 1 - B - \int_\pi^B \frac{\gamma \lambda_1 F(x, B_*)(1-x)[x/(1-x)]^\gamma}{\lambda_1 + (\lambda_0 - \lambda_1)x}\, dx,$$

$$(4.17) \quad F(\pi, B) = \frac{1}{A(\pi)\pi(1-\pi)}\left(C(\pi, B) - \int_\pi^B \frac{C(x,B)G(x)\,dx}{A(x)G(\pi)x(1-x)}\right),$$

$$(4.18) \quad A(\pi) = \frac{\lambda \lambda_0 \lambda_1 - (\lambda_0 - \lambda_1)\pi}{\pi[\lambda_1 + (\lambda_0 - \lambda_1)\pi]},$$

$$(4.19) \quad C(\pi, B) = \frac{1-B}{\gamma(\gamma-1)}\left(\frac{1-B}{B}\right)^{\gamma-1} - c(\lambda_0 - \lambda_1)\left(\frac{1-\pi}{\pi}\right)^{\gamma-1},$$

$$(4.20) \quad G(\pi) = \begin{cases} \left|\dfrac{\lambda \lambda_0 \lambda_1 - (\lambda_0 - \lambda_1)\pi}{(\lambda_0 - \lambda_1 - \lambda \lambda_0 \lambda_1)(1 - \pi)}\right|^a \dfrac{1}{1-\pi}, & \text{if } \dfrac{\lambda \lambda_0 \lambda_1}{\lambda_0 - \lambda_1} \neq 1, \\ \exp\left(\dfrac{\lambda_0 \pi}{(\lambda_1 - \lambda_0)(1 - \pi)}\right)\dfrac{1}{1-\pi}, & \text{if } \dfrac{\lambda \lambda_0 \lambda_1}{\lambda_0 - \lambda_1} = 1, \end{cases}$$

$$(4.21) \quad a = \frac{\lambda_1(1 + \lambda \lambda_0)}{\lambda_0 - \lambda_1 - \lambda \lambda_0 \lambda_1} \quad \text{if } \frac{\lambda \lambda_0 \lambda_1}{\lambda_0 - \lambda_1} \neq 1,$$

for $\pi \in (0, B]$ is a unique solution of the system (4.4) + (4.5).

Let us now assume that $\lambda_0 < \lambda_1$. In this case it follows that the function $\tilde{f}(y) = f(\pi)$ with $\pi = e^y/(1+e^y)$ solves the equation

$$(4.22) \quad \left(\frac{\lambda'(1+e^y)}{e^y} - \frac{1}{\gamma(\gamma-1)}\right)\tilde{f}'(y) - \frac{\tilde{f}(y)[\gamma(1+e^y) - 1]}{\gamma(\gamma-1)(1+e^y)}$$
$$- \frac{e^{\gamma y}}{1+e^y}\int_{-\infty}^y \frac{\tilde{f}(z)(1+e^z)}{e^{\gamma z}}\, dz = -\frac{c(\lambda_0 - \lambda_1)e^y}{1+e^y} \qquad (y < \widetilde{B})$$

and satisfies (4.10), where $\gamma = \lambda_0/(\lambda_0 - \lambda_1) < 0$, $\lambda' = \lambda(\lambda_0 - \lambda_1) < 0$ and $\widetilde{B} = \log[B/(1-B)]$. It can be easily verified that the system (4.22) + (4.10) has a unique solution which is given by

$$(4.23) \quad \tilde{f}(y) = \frac{1}{1+e^{\widetilde{B}}} + \int_{-\infty}^y \frac{\gamma(\gamma-1)\widetilde{F}(z)e^{\gamma z}}{\gamma(1+e^z) - 1}\, dz,$$

$$(4.24) \quad \widetilde{F}(y) = -\frac{c(\lambda_0 - \lambda_1)}{\widetilde{A}(y)}\left(e^{-(\gamma-1)y} + \int_{-\infty}^y \frac{e^{-(\gamma-1)z}}{\widetilde{A}(z)}\frac{\widetilde{G}(z)}{\widetilde{G}(y)}\, dz\right)$$

for $y \leq \widetilde{B}$, where $\widetilde{A}(y)$ and $\widetilde{G}(y)$ are defined in (4.13) and (4.15), respectively. Using (4.23) + (4.24) and (4.13) + (4.15) we may therefore conclude that the function $f(\pi; B) = \tilde{f}(y)$ given by (4.16) with

$$(4.25) \quad F(\pi) = -\frac{c(\lambda_0 - \lambda_1)}{A(\pi)\pi(1-\pi)}\left(\left(\frac{1-\pi}{\pi}\right)^{\gamma-1} + \int_0^\pi \frac{G(x)(1-x)^{\gamma-2}}{A(x)G(\pi)x^\gamma}\, dx\right)$$

for $\pi \in (0, B]$ is a unique solution of the system $(4.4) + (4.5)$.

Taking into account the facts proved above, we are now ready to formulate the main assertion of the section.

THEOREM 4.1. *Suppose that the observed process $X$ is given by* (4.1). *Then in the Bayesian problem of quickest disorder detection* $(2.2) + (2.3)$ *the value function $V(\pi)$ coincides with the function*

$$(4.26) \qquad V_*(\pi) = \begin{cases} f(\pi; B_*), & \pi \in (0, B_*), \\ 1 - \pi, & \pi \in [B_*, 1], \end{cases}$$

$[\text{with } V_*(0) = f(0+; B_*)]$ *and the optimal stopping time $\tau_*$ is explicitly given by* (3.5), *where $f(\pi; B)$ and the boundary $B_*$ are specified as follows:*

(i) *If $\lambda_0 > \lambda_1$ and $c > 1/\lambda_1 - 1/\lambda_0 - \lambda$, then $f(\pi; B)$ is given by* $(4.16) + (4.17)$ *and $B_* = \overline{B} \equiv \lambda/(\lambda + c)$.*

(ii) *If $\lambda_0 > \lambda_1$ and $c = 1/\lambda_1 - 1/\lambda_0 - \lambda$, then $f(\pi; B)$ is given by* $(4.16) + (4.17)$ *and $B_* = \overline{B} = \widehat{B} \equiv \lambda\lambda_0\lambda_1/(\lambda_0 - \lambda_1)$.*

(iii) *If $\lambda_0 > \lambda_1$ and $c < 1/\lambda_1 - 1/\lambda_0 - \lambda$, then $f(\pi; B)$ is given by* $(4.16) + (4.17)$ *and $B_* > \overline{B}$ is a unique root of $H(B_*) = 0$, where we set*

$$(4.27) \qquad H(B) = \int_{\widehat{B}}^{B} \frac{C(x, B)G(x)}{A(x)x(1-x)} dx.$$

(iv) *If $\lambda_0 < \lambda_1$, then $f(\pi; B) = f(\pi)$ is given by $(4.16) + (4.25)$ and $B_*$ is uniquely determined from the equation*

$$(4.28) \qquad f'(B_*) = -1.$$

PROOF. (i) and (ii) In these cases the conditions $(3.6) + (3.7)$ are satisfied and thus $\overline{B} \leq \widehat{B}$. Hence, by Lemma 3.1 we get that $B_*$ coincides with $\overline{B}$ and, by means of the uniqueness arguments for solutions of the first-order ordinary differential equations, we may conclude that $V_*(\pi) = V(\pi)$ for all $\pi \in [0, 1]$.

(iii) In this case we have $\widehat{B} < \overline{B}$, and thus, according to Remark 3.4, we see that the optimal boundary $B_*$ is located to the right of $\widehat{B}$. Taking an arbitrary $B$ from $(\widehat{B}, 1)$, by means of the arguments above we obtain that the function $f(\pi; B)$ from $(4.16) + (4.17)$ is a unique solution of the system $(4.4)$–$(4.6)$ for $\pi \in (\widehat{B}, B]$. Observe that in the given case there exists a unique point $B' \in (\widehat{B}, 1)$ such that $\lim_{\pi \downarrow \widehat{B}} f(\pi; B) = \pm \infty$ for $B \in (\widehat{B}, B') \cup (B', 1)$ and $\lim_{\pi \downarrow \widehat{B}} f(\pi; B')$ is finite. Hence $f(\pi; B)$ together with $F(\pi, B)$ from $(4.17)$ can be uniquely extended to the interval $(0, \widehat{B}]$, where by l'Hôpital's rule, we may let $F(\widehat{B}, B') = F(\widehat{B}\pm, B')$ and thus $f'(\widehat{B}; B') = f'(\widehat{B}\pm; B') \equiv -c\lambda_1^2/(\lambda_0 - \lambda_1 - \lambda\lambda_0\lambda_1)$. Then from $(4.16) + (4.17)$ it follows that $B'$ can be characterized by means of $H(B') = 0$, where $H(B)$ is defined in (4.27). Since $H(\widehat{B}+) = +0$ and the derivative $H'(B) > 0$ for $B \in (\widehat{B}, \overline{B})$ and $H'(B) < 0$ for $B \in (\overline{B}, 1)$, the function $H(B)$ increases on $(\widehat{B}, \overline{B})$ and decreases on $(\overline{B}, 1)$. Thus, by virtue of the property $\lim_{B \uparrow \infty} H(B) = -\infty$, we get that $B'$ belongs to the interval $(\overline{B}, 1)$ and $H(B') = 0$ has a unique solution.

Summarizing the facts proved above, we see that the value function $V(\pi)$ and the optimal boundary $B_*$ should necessarily solve the system (4.4)–(4.6) and there is only one point $B'$ such that the solution $f(\pi; B')$ taken at $\pi = \widehat{B}$ is finite. We may therefore conclude that $B_*$ coincides with $B'$ and the uniqueness argument for solutions of first-order differential equations implies that $V_*(\pi) = V(\pi)$ for all $\pi \in [0, 1]$, thus proving the claim.

(iv) Taking into account the fact that in this case the process $(\pi_t)$ can increase only continuously, following the arguments in [17], Chapter IV.4, and [10] we may guess that the smooth-fit condition (4.7) is satisfied and thus (4.28) holds. Using straightforward calculations it is shown that $f''(\pi) < 0$ for $\pi \in (0, 1)$; hence, the function $f(\pi)$ from (4.16) + (4.25) is concave on $[0, 1]$ and its derivative $f'(\pi)$ is decreasing on $(0, 1)$. Therefore, by virtue of the facts that $f'(0+) = 0$ and $f'(1-) = -\infty$, we may conclude that (4.28) admits a unique solution.

Let us now show that the function $V_*(\pi)$ defined in (4.26) + (4.16) + (4.25) coincides with the value function $V(\pi)$ and that $B_*$ being a unique root of (4.28) is an optimal stopping boundary. For this, applying Itô's formula, we get

$$(4.29) \qquad V_*(\pi_t) = V_*(\pi) + \int_0^t (\mathbb{L} V_*)(\pi_{s-}) \, ds + M_t^*,$$

where the process $(M_t^*)$ defined by

$$(4.30) \quad M_t^* = \int_0^t \int_0^\infty \left[ V_* \left( \frac{\pi_{s-} \exp(-\lambda_1 x)}{\pi_{s-} \exp(-\lambda_1 x) + (1 - \pi_{s-}) \exp(-\lambda_0 x)} \right) - V_*(\pi_{s-}) \right]$$
$$\times (\mu^X(ds, dx) - (\pi_{s-} \exp(-\lambda_1 x) + (1 - \pi_{s-}) \exp(-\lambda_0 x)) \, ds \, dx)$$

is a martingale under $P_\pi$ with respect to $\mathbf{F}^X$.

Since $V_*(\pi)$ is a bounded function, from (4.30) by means of the optional sampling theorem we get that $E_\pi[M_\tau^*] = 0$ for all $\tau$ from $\mathcal{M}(\pi)$. Thus, taking the expectation on both sides in (4.29) with $\tau$ instead of $t$ and using the fact that a direct verification yields $(\mathbb{L} V_*)(\pi) \geq -c\pi$ and $V_*(\pi) \leq 1 - \pi$, we obtain

$$(4.31) \qquad V_*(\pi) \leq E_\pi \left[ 1 - \pi_\tau + c \int_0^\tau \pi_t \, dt \right]$$

for all $\tau$ from the class $\mathcal{M}(\pi)$, and hence $V_*(\pi) \leq V(\pi)$ for all $\pi \in [0, 1]$.

Observe that straightforward calculations above imply that the function $V_*(\pi)$ and the boundary $B_*$ solve the system (4.4)–(4.6); hence we have $V_*(\pi_{\tau_*}) = 1 - \pi_{\tau_*}$ and $(\mathbb{L}V_*)(\pi_t) = -c\pi_t$ for all $0 \le t \le \tau_*$. Therefore, taking the expectation on both sides in (4.29) with $t$ replaced by $\tau_*$ and using the obvious fact that $\tau_*$ belongs to $\mathcal{M}(\pi)$, we see that the equality in (4.31) is attained at $\tau = \tau_*$. This implies that $V_*(\pi) = V(\pi)$ for all $\pi \in [0,1]$ and that $B_*$ is an optimal stopping boundary. Thus the proof is complete. $\square$

REMARK 4.2. We observe that in case (i) of Theorem 4.1 we can verify that $f'(B_*-; B_*) = -1$ and in the case (iv) we have proved that (4.28) holds, so that the smooth-fit condition (4.7) is satisfied. This can be explained by the facts that the process $(\pi_t)$ may pass through $B_*$ continuously and that (4.4) has no singularity point. On the other hand, in case (ii) it is shown that $f'(B_*-; B_*) = -c\lambda_1^2/(\lambda_0 - \lambda_1 - \lambda\lambda_0\lambda_1) > -1$ and in case (iii) it can be also proved that the smooth-fit condition (4.7) breaks down. This can be explained by means of the facts that the process $(\pi_t)$ may pass through $B_*$ for the first time only by jumping and that (4.4) has a singularity point $\widehat{B}$.

REMARK 4.3. We note that the function $f(\pi; B)$ for different $B \in (0,1)$ and the function $V_*(\pi)$ in cases (i)–(iv) look the same as in [10], Figures 2–5.

**5. Solution of the variational problem for a compound Poisson process with exponential jumps.** Let us first note that if $\alpha \ge 1 - \pi$, then letting $\hat{\tau} = 0$ we get $P_\pi[\hat{\tau} < \theta] = P_\pi[\theta > 0] = 1 - \pi \le \alpha$ and $E_\pi[\hat{\tau} - \theta]^+ = 0$, whence it is seen that $\hat{\tau} = 0$ is optimal in the formulation $(2.4) + (2.5)$.

Assuming that $0 < \alpha < 1 - \pi$ and following the arguments from [17], pages 198–200, we further show that the solution of the variational problem $(2.4) + (2.5)$ can be obtained using the solution of the Bayesian problem. For this, let us introduce the function

(5.1) $\quad u(\pi; B_*) = P_\pi[\tau_* < \theta] \qquad (= E_\pi[1 - \pi_{\tau_*}]).$

To find an explicit expression for the function $u(\pi; B)$ in the case when $\lambda_0 > \lambda_1$, we observe that, by virtue of the strong Markov property, it should solve the system

(5.2) $\quad (\mathbb{L}u)(\pi; B) = 0 \qquad (0 < \pi < B),$

(5.3) $\quad u(\pi; B) = 1 - \pi \qquad (B \le \pi \le 1).$

By means of the same arguments as in the text that accompanies the formulas (4.9)–(4.21), it is shown that the system $(5.2) + (5.3)$ admits the unique solution

(5.4) $\quad u(\pi; B) = 1 - B - \int_\pi^B \frac{\gamma \lambda_1 D(x, B)(1-x)}{\lambda_1 + (\lambda_0 - \lambda_1)x} \left(\frac{x}{1-x}\right)^\gamma dx,$

(5.5) $\quad D(\pi, B) = \frac{1-B}{\gamma(\gamma-1)A(\pi)\pi(1-\pi)} \frac{G(B)}{G(\pi)} \left(\frac{1-B}{B}\right)^\gamma$

for $\pi \in (0, B)$, $\pi \neq \widehat{B}$, where $\gamma = \lambda_0/(\lambda_0 - \lambda_1) > 1$, the functions $A(\pi)$ and $G(\pi)$ are given by (4.18) and (4.20), respectively, and by l'Hôpital's rule, we can let $D(\widehat{B}, B) = D(\widehat{B}\pm, B) \equiv 0$ as well as $u(0; B) = u(0+; B)$.

It is not difficult to verify that $\partial u(\pi; B)/(\partial B) < 0$ for $B \in (\pi, 1)$, so that the function $u(\pi; B)$ is strictly decreasing on $(\pi, 1)$ for $0 < \pi < 1 - \alpha$ fixed. Therefore, by virtue of the obvious facts that $u(\pi; 0) = 1 - \pi$ and $u(\pi; 1) = 0$, we may conclude that there exists a point $B(\alpha) \leq 1 - \alpha$ that is a unique solution of the equation

$$u(\pi; B(\alpha)) = \alpha. \tag{5.6}$$

Let us now formulate the main result of the section.

THEOREM 5.1. *Suppose that the observed process $X$ is given by (4.1). Then in the variational problem of quickest disorder detection (2.4) + (2.5), the optimal stopping time $\hat{\tau}$ is explicitly given by*

$$\hat{\tau} = \inf\{t \geq 0 | \pi_t \geq B(\alpha)\}, \tag{5.7}$$

*where the boundary $B(\alpha) \leq 1 - \alpha$ is specified as follows:*

   (i) *If $0 < \alpha < 1 - \pi$ and $\lambda_0 > \lambda_1$, then $B(\alpha)$ is a unique root of (5.6).*
   (ii) *If $\alpha \geq 1 - \pi$ or $\lambda_0 < \lambda_1$, then $B(\alpha) = 1 - \alpha$.*

PROOF. (i) Let us consider the function $B_* = B_*(c)$ as an optimal boundary in the corresponding Bayesian problem which is uniquely determined from parts (i)–(iii) of Theorem 4.1. It can be easily shown that $B_*(c)$ is continuous and strictly decreasing on $(0, \infty)$, and it satisfies $\lim_{c \downarrow 0} B_*(c) = 1$ and $\lim_{c \uparrow \infty} B_*(c) = 0$. Then there exists a constant $c(\alpha)$ such that $B(\alpha) = B_*(c(\alpha))$ and by the definition (2.2), we have

$$P_\pi[\hat{\tau} < \theta] + c(\alpha) E_\pi[\hat{\tau} - \theta]^+ \leq P_\pi[\tau < \theta] + c(\alpha) E_\pi[\tau - \theta]^+ \tag{5.8}$$

for all stopping times $\tau$. Since from (5.6) together with (5.1) and (3.5) it is seen that $P_\pi[\hat{\tau} < \theta] = \alpha$, we may thus conclude that (5.8) directly yields

$$c(\alpha) E_\pi[\hat{\tau} - \theta]^+ \leq c(\alpha) E_\pi[\tau - \theta]^+ \tag{5.9}$$

for all $\tau$ from $\mathcal{M}(\pi, \alpha)$. Therefore, by virtue of the obvious fact that $c(\alpha) > 0$ for $0 < \alpha < 1 - \pi$, we obtain that $\hat{\tau}$ from (5.7) is optimal in (2.5).

(ii) Since whenever $\lambda_0 < \lambda_1$, the process $(\pi_t)$ can increase only continuously, we get that $\{\pi_{\hat{\tau}} \geq B(\alpha)\} = \{\pi_{\hat{\tau}} = B(\alpha)\}$, and from (5.1) it thus follows that in this case we have $u(\pi; B) = 1 - B$. Hence, from (5.6) it is seen that $B(\alpha) = 1 - \alpha$, and the arguments from the previous part (i) complete the proof. $\square$

ON THE DISORDER PROBLEM 13**Acknowledgments.** I am grateful to A. N. Shiryaev and G. Peskir for the statement of the problem and for many helpful discussions. I am thankful to the Editor for the encouragement to prepare the revised version, and am obliged to an Associate Editor and a referee for many useful suggestions which are incorporated into the final version of the paper.## REFERENCES

[1] BARNDORFF-NIELSEN, O. E. (1995). Normal inverse Gaussian processes and the modelling of stock returns. Research Report 300, Dept. Theoretical Statistics, Aarhus Univ.
[2] BERRY, D. A. and FRISTEDT, B. (1985). *Bandit Problems*: *Sequential Allocation of Experiments*. Chapman and Hall, London. MR813698
[3] CARLSTEIN, E., MÜLLER, H.-G. and SIEGMUND, D., eds. (1994). *Change-Point Problems*. IMS, Hayward, CA. MR1477909
[4] DAVIS, M. H. A. (1976). A note on the Poisson disorder problem. *Banach Center Publ.* **1** 65–72.
[5] DYNKIN, E. B. (1963). The optimum choice of the instant for stopping a Markov process. *Soviet Math. Dokl.* **4** 627–629. MR193670
[6] GAL'CHUK, L. I. and ROZOVSKII, B. L. (1971). The "disorder" problem for a Poisson process. *Theory Probab. Appl.* **16** 712–716. MR297028
[7] JACOD, J. and SHIRYAEV, A. N. (1987). *Limit Theorems for Stochastic Processes*. Springer, Berlin. MR959133
[8] KOLMOGOROV, A. N., PROKHOROV, YU. V. and SHIRYAEV, A. N. (1990). Methods of detecting spontaneously occurring effects. *Proc. Steklov Inst. Math.* **1** 1–21.
[9] MORDECKI, E. (1999). Optimal stopping for a diffusion with jumps. *Finance Stochastics* **3** 227–236. MR1810074
[10] PESKIR, G. and SHIRYAEV, A. N. (2002). Solving the Poisson disorder problem. In *Advances in Finance and Stochastics. Essays in Honour of Dieter Sondermann* (K. Sandmann and P. Schönbucher, eds.) 295–312. Springer, Berlin. MR1929384
[11] SATO, K. I. (1999). *Lévy Processes and Infinitely Divisible Distributions*. Cambridge Univ. Press. MR1739520
[12] SHIRYAEV, A. N. (1961). The detection of spontaneous effects. *Soviet Math. Dokl.* **2** 740–743.
[13] SHIRYAEV, A. N. (1961). The problem of the most rapid detection of a disturbance in a stationary process. *Soviet Math. Dokl.* **2** 795–799.
[14] SHIRYAEV, A. N. (1963). On optimum methods in quickest detection problems. *Theory Probab. Appl.* **8** 22–46.
[15] SHIRYAEV, A. N. (1965). Some exact formulas in a "disorder" problem. *Theory Probab. Appl.* **10** 348–354.
[16] SHIRYAEV, A. N. (1967). Two problems of sequential analysis. *Cybernetics* **3** 63–69. MR272119
[17] SHIRYAEV, A. N. (1978). *Optimal Stopping Rules.* Springer, Berlin. MR468067
[18] SHIRYAEV, A. N. (1999). *Essentials of Stochastic Finance.* World Scientific, Singapore. MR1695318

Russian Academy of Sciences
Institute of Control Sciences
Profsoyuznaya Street 65
117997 Moscow
Russia
e-mail: gapeev@cniica.ru